\documentclass{article}

\usepackage{amssymb,dsfont,latexsym}

\newtheorem{thm}{Theorem}

\newtheorem{lem}[thm]{Lemma}

\newcommand\enu[1]{\smallskip\newline\makebox[5mm][l]{\rm(#1)}}

\newcommand\bp{\noindent{\it Proof.}\ }

\newcommand\Tr{{\rm Tr}}

\begin{document}

\author{Erling St{\o}rmer}

\date{27-2-2009 }

\title{Decomposable and atomic projection maps}

\maketitle
\begin{abstract}
It is shown that a trace invariant projection map, i.e. a positive unital idempotent map, 
of a finite dimensional $C^*-$algebra into itself is nondecomposable if and only if it
is atomic, or equivalently not the sum of a 2-positive and a 2-copositive map.  In
particular projections onto spin factors of dimension greater than 6 are atomic.

\end{abstract}

\section*{Introduction}
The classification theory for positive linear maps between $C^*-$algebras has been
progressing slowly since Stinepring \cite{St} introduced the class of completely
positive maps in 1955.  Even in the simple case when both $C^*-$algebras are the
complex $3\times 3$ matrices $M_{3}(\0C)$, the classification problem is still open.  
So far main emphasis has been on completely positive maps, copositive maps, i.e. those
which are the composition of a completely positive map and an anti-automorphism, or
more generally,  k-positive maps, which are maps $\phi$ such that $\phi\otimes \iota_{k}$
is positive, where $\iota_{k}$ is the identity on $M_{k}(\0C)$, and the corresponding 
k-copositive maps. These maps are reasonably well understood, see e.g. \cite{SSZ}.
Furthermore maps which are sums of completely positive maps and copositive maps,
called \textit{decomposable maps}, and those which are not, called \textit
{nondecomposable maps}, have attracted much attention.

A more general class of which we have a reasonably good understanding is the cone of
\textit{(k,m)-decomposable maps}, which are those which are the sum of a k-positive and an
m-copositive map, including the case when one of them is the zero map.  They are all 
contained in the cone of (2,2)-decomposable maps.  It is thus natural to consider maps  
which are not (2,2)-decomposable, called \textit{atomic} in the literature.  It was shown
by Ha \cite{H} that the natural generalizations of the Choi map \cite{Cho} are atomic.
In the present paper we shall show that trace preserving projection maps are atomic if and only if  they
are nondecomposable.  In \cite {S2} it was shown that the projetion maps onto spin 
factors of dimension greater than 6 are nondecomposable, hence we obtain atomic maps
in arbitrary large dimensions.  We shall at the end apply the above results to give a sufficient
condition for a positive map to be be atomic, in terms of algebraic properties of its fixed
point set.

\section {Atomic projection maps}
If $B$ is a unital $C^*-$algebra a positive linear map $P\colon B \to B$ is called a 
\textit {projection map} if $P(1) = 1$ and $P^2 = P$.  The theory of such maps is 
intimately related to the theory of Jordan algebras of self-adjoint operators under the 
Jordan product $a\circ b = \frac{1}{2}(ab + ba)$ for $a,b\in B_{sa}$, the self-adjoint
operators in $B$.  See \cite{HS} for the theory of Jordan algebras of self-adjoint 
operators.  In fact, if $P$ is faithful, so $P(a)\neq0$ for $a>0$, then $A = P(B_{sa})$ is a 
Jordan subalgebra of $B_{sa},$ \cite{ES}.  Furthermore, by \cite{S2} $P$ is 
decomposable if and only if $A$ is a \textit{reversible Jordan algebra}, i.e. $A$ is closed 
under symmetric products $a_{1}a_{2}...a_{n} + a_{n}a_{n-1}...a_{1},$ with $a_{i}\in A.$
In the present section we shall sharpen the above result by replacing decomposable by
(2,2)-decomposable maps.  For simplicity we shall assume the $C^*-$algebra $B$ is
acting on a finite dimensional Hilbert space $H$ and denote the usual trace on $B(H)$ by
$Tr$.  

\begin{thm}\label{thm 1}
Let $B$ be a $C^*-$algebra acting on a finite dimensional Hilbert space.  Let $P\colon
B\to  B$ be a trace preserving unital projection map. Then we have:
\enu{i} $P$ is (2,2)-decomposable if and only if $P$ is decomposable.
\enu{ii} $P$ is atomic if and only if $A= P(B_{sa})$ is a nonreversible $Jordan $algebra.
\end{thm}

The proof will be divided into some lemmas.  We first recall that a \textit{spin system}
in $B(H)$ is a set $(s_{i})_{i\in I}$ of symmetries, i.e. self-adjoint unitaries in $B(H),$
such that $s_{i}\circ s_{j} = 0$ for $i \neq j.$  The real linear span of $1$ and the $s_{i}$ 
is a Jordan algebra called a \textit{spin factor}, see \cite {HS} Chapter 6. If $H$ is 
finite dimensional there is a canonical positive trace preserving projection map $P$ of $B(H)$
onto $A+ iA,$ see \cite{ES}.  For simplicity we often write $P\colon B(H)\to A.$

\begin{lem}\label{lem 2}
Let $A\subset B(H)$ be a spin factor.  Let $e_{1},...,e_{k}$ be nonzero minimal projections
in the commutant of $A$ with sum $1$. Let 
$$
P_{i}\colon  e_{i}B(H)e_{i} \to Ae_{i}
$$
be the canonical projection map, and let $P\colon B(H) \to A$ be the canonical projection map. 
Let $\alpha_{i}\colon Ae_{i} \to A$ by $\alpha(ae_{i}) = a.$  Then $\alpha_{i}$ is an 
isomorphism, and
$$
P(a) = \sum    \Tr(e_{i})\alpha_{i}(P_{i}(e_{i}a e_{i})).
$$
Hence $P$ is decomposable if and only if each $P_{i}$ is decomposable.
\end{lem}
\bp
Since $A$ is a simple Jordan algebra $\alpha_{i}$ is an isomorphism. Let $(s_{j})$ be a spin 
system generating $A.$ Then $(s_{j}e_{i})$ is a spin system generating $Ae_{i}$. Let
$c_{i} = \Tr(e_{i})^{-1}.$  Then $c_{i}\Tr$ is the tracial state on $e_{i}B(H)e_{i}$, and
the canonical projection $P_{i}$ being the orthogonal projection of $e_{i}B(H)e_{i}$ with 
respect to the Hilbert Schmidt structure onto the subspace generated by $Ae_{i}$ is given by 
$$
P_{i}(e_{i}ae_{i}) = c_{i}\Tr(e_{i}ae_{i})e_{i} + c_{i} \sum_{j} \Tr(e_{i}ae_{i}s_{j}e_{i})s_{j}e_{i}
                             = c_{i}P(e_{i}ae_{i})e_{i}.
 $$
 Thus we have
 $$
 \alpha_{i}(P_{i}(e_{i}ae_{i})) = \alpha_{i}(c_{i}P(e_{i}ae_{i})e_{i}) = c_{i}P(e_{i}ae_{i}).
 $$
 Let $E\colon B(H) \to \bigoplus e_{i}B(H)e_{i}$, given by $E(x) = \sum_{i} e_{i}xe_{i}$.  Then $E$ is completely positive, 
and by uniqueness of $P$ as the trace invariant projection of $B(H)$ onto $A$,  $P$ is the restriction
$$
P =( P\vert \bigoplus e_{i}B(H)e_{i}) \circ E.
$$
Thus, if $a\in B(H)$ then
$$
P(a) = \sum_{i} P(e_{i}ae_{i}) = \sum_{i} c^{-i} \alpha_{i}P_{i}(e_{i}ae_{i})
= \sum_{i} \Tr(e_{i})\alpha_{i}(P_{i}(e_{i}ae_{i}))
$$
 Since $E$ and $\alpha_{i}$ are completely positive, $P$ is decomposable if and only if each $P_{i}$ is
 decomposable.  The proof is complete.
\medskip
 
 In the case when $A$ is a spin factor we have reduced the proof of the theorem to the case when $A$
 is irreducible, $C^*(A) = B(H)$, or equivalently the commutant of $A$ is the scalars.

\begin{lem}\label{lem3}
Let $R_{2}$ denote the real symmetric $2\times 2$ matrices.  Let $\phi\colon M_{2}(\0C) \to
M_{2}(\0C)$ be a positive map such that if $\iota_{R_{2}}$  denotes the identity map on $R_{2},$ then 
$\iota_{R_{2}} \geq \mu \phi$ for some $\mu \geq 0.$  Then $\phi = \lambda  \iota_{R_{2}}$ for some 
$\lambda \geq 0.$
\end{lem}
 \bp
 Let $e$ be a 1-dimensional projection in $R_{2}.$  Then $e = \iota_{R_{2}}(e) \geq \mu \phi(e)\geq 0,$ so 
 $\phi(e)= \alpha e, \alpha \geq 0.$ Similarly $\phi(1-e)= \beta (1-e), $  since $1-e$ is also 
 1-dimensional.  Thus $\phi(1)$ belongs to the maximal abelian subalgebra of $R_{2}$ generated by $e.$
 This holds for all minimal projections $e$ in $R_{2},$ which is possible only if $\phi(1)=\lambda 1$
 for some $\lambda \geq 0$. Since $\lambda 1 = \alpha e +  \beta (1-e), \alpha = \beta = \lambda,$
 i.e. $\phi(e)=\lambda e$ for all minimal projections $e\in R_{2},$ so that $\phi = \lambda  \iota_{R_{2}}.$
 The proof is complete.
 
 \begin{lem}\label{lem4}
 Let $A\subset B(H)$ be a spin factor and $P$ the canonical projection of $B(H)$ onto $A.$  Suppose
 $\phi\colon B(H)\to B(H)$ is a positive linear map such that $\phi \leq \mu P$ for some $\mu > 0.$
 Then the restriction $\phi\vert A = \lambda \iota_{A}$  for some $\lambda \geq 0,$ where $\iota_{A}$
 is the identity map on $A.$
 \end{lem}
  \bp
  Let $s_{i}$ and $s_{j}, i \neq j,$  belong to the spin system generating $A,$ and let $F_{ij}$ denote the canonical trace 
  preserving conditional expectation of $B(H)$ onto the $C^*-$algebra $C^*(s_{i},s_{j})$ generated by
  $s_{i}$ and $s_{j}.$ Since $span\{1,s_{i},s_{j}\} \simeq R_{2}$, by Lemma 3 the restriction 
  $F_{ij}\circ \phi\vert_{R_{2}} =\lambda \iota_{R_{2}}$  for some $\lambda \geq 0.$ Thus we have
  $F_{ij}\circ\phi(1)=\lambda 1,$ and $F_{ij}\circ \phi(s_{i})=\lambda  s_{i}$ and similarly for $s_{j}.$ Since this holds
  for all $j\neq i,$ it follows that $F_{ij}\circ\phi(1)=\lambda 1,$ and $F_{ij}\circ\phi(s_{i}) = \lambda s_{i}$
  for all $i\neq j$. Thus $\phi(1)=\lambda 1 + a_{ij}$ with $a_{ij}$ orthogonal to $C^*(s_{i},s_{j})$ in the
  Hilbert Schmidt structure.  But then  $a_{ij}=a_{kl} = a$ for all $i,j,k,l,$  and so $a$ is orthogonal to 
  $A^2 = \{xy: x,y \in A\}.$  Similarly $\phi(s_{i})=\lambda s_{i} + b$ with $b$ orthogonal to $A^2.$
  Scaling $\phi$ we may assume 
  $$
  \phi(1)= 1 + a,   \phi(s_{i}) = s_{i} +b_{i},  a, b_{i} \perp A^2.
  $$
   Let $s= s_{i}$ and $e=e_{i}=\frac{1}{2}(1 + s_{i})=\frac{1}{2}(1+ s).$  Then $e$ is a projection in $A.$  Since $\phi\leq\mu P,$ 
  $\phi(e)\leq \mu e.$  Since $f=\frac{1}{2}(1 - s) =1 - e$ is a projection 
  orthogonal to  $e$, we have $\phi(e)f=0,$ and similarly $\phi(f)e = 0.$ Calculating we get, using that $s^2 = 1$
  $$
  0 = (\phi(1+ s ) )(1 - s) = ((1 + a )+ (s + b))(1 - s) =  a + b - as -bs,
  $$
  $$
  0 = (\phi(1 - s))(1 + s)  = ((1 + a ) - (s + b))(1 + s)=  a - b +as -bs.
  $$
Adding these two equations we get $0 = a - bs,$ hence $\phi(1) = 1  + bs,$  and $ \phi(s)= s  + b = 
(1 + bs)s = \phi(1)s.$ Thus the product of the two self-adjoint operators $\phi(1)$ and $s$ is self-adjoint, hence they
commute.  Since $s = s_{i}$ was an arbitrary symmetry in the spin system spanning $A, \phi(1)$ belongs to the 
commutant of $A.$  But $A$ was assumed to be irreducible, so  $bs = 0 = b$, so that 
$\phi(1)= 1$, and $\phi(s_{i})= s_{i}$ for all i.  This completes the proof of the lemma.

\begin{lem}\label{lem5}
Let $A$ be an irreducible spin factor acting on the finite dimensional Hilbert space $H$.  Let $\phi, \psi\colon 
B(H)\to B(H)$ be unital maps with $\phi$ 2-positive and $\psi$ 2-copositive such that their restrictions to $A$
are the identity map.  Let $s,t$ be distinct symmetries in the spin system spanning $A$. Then
$$
\phi(st) = st,   \psi(st) = ts.
$$
\end{lem}
\bp
Let $x=s + it$.  Since $st = -ts$ we have 
$$
x^*x = 2(1 + ist),   xx^* = 2(1 - ist).
$$
Hence $x^*x + xx^* = 4\cdot 1\in A.$  Since $\phi$ is 2-positive it satisfies the Schwarz inequality \cite{Cho}, Cor.
2.8, hence $\phi(x^*x)\geq \phi(x)^*\phi(x) = x^*x$, and $\phi(xx^*)\geq \phi(x)\phi(x)^* = xx^*.$ Thus
$$
0 = \phi(x^*x + xx^*) - 4\cdot 1 = \phi(x^*x) + \phi(xx^*) - 4\cdot 1 \geq x^*x + xx^* - 4\cdot 1 =0,
$$
Hence $\phi(x^*x) = x^*x$, and $\phi(xx^*) = xx^*.$  In particular  $\phi(st) = st.$ 

Since $\psi $ is 2-copositive , $\psi(x^*x)\geq \psi(x)\psi(x)^* = xx^*,$ so by the above argument applied
to $\psi$ we get $\psi(st) = ts,$  completing the proof.

\begin{lem}\label{lem6}
Let $\phi,\psi\colon B(H)\to B(H)$ be unital maps with $\phi$ 2-positive and $\psi$ 2-copositive.
Let $a\in B(H).$  Then we have:
\enu{i} If $\phi(aa^*) = \phi(a)\phi(a)^*,$ then $\phi(ab)= \phi(a)\phi(b) \ \forall  b\in B(H).$
\enu{ii} If $\psi(a^*a) = \psi(a)\psi(a)^*,$ then $\psi(ba)= \psi(a)\psi(b)\ \forall b\in B(H).$
\end{lem}
\bp
Let 
$$
<x,y> = \phi(x y^*) - \phi(x) \phi(y)^*,   x,y \in B(H).
$$
Then $< , >$ is an operator valued sesquilinear form such that for all states $\omega$ on $B(H), <x,y>_{\omega}=
\omega(<x,y>)$ is a sesquilinear form.  Thus by the Cauchy - Schwarz inequality
$$
\vert<x,y>_{\omega}\vert
 \leq \omega(\phi(xx^*) - \phi(x)\phi(x)^*)^{\frac{1}{2}}\omega(\phi(yy^*) - \phi(y)\phi(y)^*)^{\frac{1}{2}}.
$$
Hence if $a$ is as in the statement of (i), then $<a,b>_{\omega} = 0$ for all $\omega$, hence $<a,b> =0,$ i.e.
$\phi(ab^*) = \phi(a)\phi(b^*)$  for all $b\in B(H)$, proving (i).

(ii) In this case we consider the operator valued sesquilinear form
$$
<x,y> = \psi(x y^*) - \psi(y)^*\psi(y),
$$
and we use the same argumens as in (i).  The proof is complete.

\medskip
Let $V_{k}$ denote the spin factor generated by a spin system consisting of k symmetries, defined as in
\cite{HS}, section 6.2.

\begin{lem}\label{lem7}
Let $P$ be the canonical projection of $B(H)$ onto the spin factor $A.$ Suppose $P$ is (2,2)-decomposable.
Then $P$ is decomposable, and $A$ is one of the spin factors $V_{2},V_{3}, V_{5}$.
\end{lem}
\bp
By Lemma 2 we may assume the $C^*-$algebra generated by $A, C^*(A)= B(H).$ Suppose $P = \phi + \psi$
with $\phi$ 2-positive and $\psi$ 2-copositive. By Lemma 4 we can replace $\phi$ and $\psi$ by maps 
which are the identity map on $A,$ and asssume $P = \lambda \phi + (1 - \lambda)\psi,  0\leq \lambda \leq 1.$
Then by Lemma 5 if $s_{k} \neq s_{j}$ are symmetries in the spin system spanning $A$ then $\phi(s_{k}s_{j})=
s_{k}s_{j},$ and $\psi(s_{k}s_{j}) = s_{j}s_{k}.$  From the proof of Lemma 5 if $x = s_{k} + is_{j},$ then
$\phi(xx^*)=xx^*,$ and $\psi(xx^*)= x^*x$, so by Lemma 6
$$
\phi(xx^*b) = \phi(xx^*)\phi(b) = xx^*\phi(b) \ \    \forall  b\in B(H).
$$
Now the monomials $s_{i_{1}}s_{i_{2}}...s_{i_{n}}$   span $C^*(A) = B(H)$ linearly.  For such monomials we 
have, since the above equation holds in particular for $s_{k}s_{j}$,
\begin{eqnarray*}
\phi(s_{i_{1}}s_{i_{2}}...s_{i_{n}}) &=& s_{i_{1}}s_{i_{2}}\phi(s_{i_{3}}...s_{i_{n}})\\
&=& s_{i_{1}}s_{i_{2}}s_{i_{3}}s_{i_{4}}\phi(s_{i_{5}}...s_{i_{n}})\\
&=& ...\\
&=& s_{i_{1}}s_{i_{2}}...s_{i_{n}}.
\end{eqnarray*}
Thus $\phi$ is the identity map.  Similarly we have
$$
\psi(s_{i_{1}}s_{i_{2}}...s_{i_{n}}) = s_{i_{n}}s_{i_{n - 1}}...s_{i_{1}},
$$
so that $\psi$ is an anti-automorphism of order 2, which is the identity on $A.$  In particular $\phi$ is
completely positive, and $\psi$ is copositive.  It follows that $P$ is decomposable.  It then follows
from \cite{S2}, Corollary 7.3, that $A$ is reversible, hence by \cite{HS}, Theorem 6.2.5, that $A$
is one of the spin factors $V_{2}, V_{3},V_{5}$.  The proof is complete. 

\medskip
\textit{Proof of Theorem 1}. Let $B$ be a $C^*-$algebra acting on the finite dimensional Hilbert 
space $H,$  and let $P$ be a unital trace preserving positive projection map of $B$ into itself. Let $A=
P(B_{sa})$.  By \cite{ES} $A$ is a Jordan subalgebra of $B_{sa}.$ Composing $P$ with the trace 
invariant conditional expectation $E$ of $B(H)$ onto $B$ we may assume that $P$ is a projection of $B(H)$
onto $A.$   This does not alter the conclusion of the theorem since $E$ is completely positive, see e.g.
\cite{str},Proposition 9.3.  Let $Z$ denote the center of $A,$ see \cite{HS}, 2.5.1, which is the
self-adjoint part of an abelian $C^*-$algebra.  Let $p$ be a minimal projection in $Z.$  Then by \cite{HS},
Proposition 5.2.17, the center of $pA = pAp = Zp =\0R p,$ so $pAp$ is a Jordan factor, also called a 
JW-factor.  It follows that $A = \bigoplus_{j } Ap_{j},$ with $p_{j}$ a minimal central projection
in $A,$ and $Ap_{j}$ is a Jordan factor.

By \cite{HS},section 6.3, a Jordan factor $C$ is either a spin factor or is reversible. In the latter case
a  projection onto $C$ is necessarily decomposable,\cite{S2}, Corollary 7.4, hence it remains to 
consider the case when $Ap_{j}$ is a spin factor.  But then Lemma 7 implies that the projection is 
decomposable if and only if it is (2,2)-decomposable, which proves part (i) of the theorem.  Part (ii)
is immediate from part (i) by again applying \cite{S2}. The proof is complete.

\medskip
A closer look at the proof shows that we have proved a slightly more general result.  Instead of assuming 
the projection map $P$ is the sum of a 2-positive map $\phi$ and a 2-copositive map $\psi$, we could
have assumed $\phi$ to satisfy the Schwarz inequality $\phi(x^*x)\geq \phi(x)^*\phi(x)$ and $\psi$
the inequality $\psi(x^*x)\geq\psi(x)\psi(x)^*$.

Theorem 1 has a natural generalization to positive maps.

If $B$ is a finite dimensional $C^*-$algebra and $\phi\colon B\to B$ is a positive unital map, then
the fixed point set $B_{\phi} = \{a\in B: \phi(a) = a\}$ has a natural structure as a Jordan algebra,
see \cite{ES}.  Furthermore, if there exists a faithful $\phi-$invariant state on $B,$ then there
exists a faithful $\phi-$invariant positive projection map $P_{\phi}\colon B\to B_{\phi}$ 
making $(B_{\phi})_{sa}$ a Jordan subalgebra of $B_{sa. }$ We then have the following extension 
of Theorem 1.

\begin{thm}\label{thm8}
Let $B$ be $C^*-$algebra acting on a finite dimensional Hilbert space.  Suppose $\phi\colon B\to B$
is a trace preserving unital positive map.  Let $P_{\phi}\colon B \to B_{\phi}$ be the 
$\phi-$invariant positive projection of $B$ onto $B_{\phi}.$  If $\phi$ is (2,2)-decomposable then 
$P_{\phi}$ is decomposable, and $(B_{\phi})_{sa}$ is a reversible Jordan algebra.  In particular 
if $(B_{\phi})_{sa}$ is nonreversible, then $\phi$ is atomic.
\end{thm}
\bp
The projection map $P_{\phi}$ is a weak limit of averages $\phi_{n} = \frac{1}{n}\sum_{0}^{n-1}\phi^{i}$.
If $\phi$ is (2,2)-decomposable, so is $\phi_{n}$, say $\phi_{n}=\alpha_{n}+\beta_{n}$ with $\alpha_{n}$
2-positive, and $\beta_{n}$ 2-copositive.  Now a subnet of $(\phi_{n})$ converges weakly to $P_{\phi}$.  If
$\alpha$ and $\beta$ are corresponding weak limit points of $\alpha_{n}$ and $\beta_{n}$ then $\alpha$
is 2-positive, and $\beta$ is 2-copositive.  Thus $P_{\phi}$  is (2,2)-decomposable, so by Theorem 1 $P_{\phi}$
is decomposable, because $\phi$ was assumed to be trace preserving, hence so is $P_{\phi},$ and therefore
$P_{\phi}$ is faithful.  By \cite{S2}, Corollary 7.3, $(B_{\phi})_{sa}$ is reversible. The proof is complete.

Department of Mathematics, University of Oslo, 0316 Oslo, Norway.

e-mail: erlings@math.uio.no

\end{document}